\documentclass[a4paper,12pt]{article}
\usepackage{amsfonts,graphics,amsmath}
\usepackage{rotating}

\newcommand{\eqnwithbr}[2]{%
\refstepcounter{equation}
\begin{trivlist}
\item[]#1 \hfill $\displaystyle #2$ \hfill (\theequation)
\end{trivlist}}

\newcommand{\eqnlabel}[2]{%
\begin{trivlist}
\item[]#1 \hfill $\displaystyle #2$ \hfill~
\end{trivlist}}

\newcommand{\mb}{\mathbf}

\renewcommand{\theenumi}{\roman{enumi}}

\makeatletter \renewcommand\p@enumi{(\theenumi)} \makeatother

\def\complaint#1{}
\def\withcomplaints{
\newcounter{mycomplaints}
\def\complaint##1{\refstepcounter{mycomplaints}%
\ifhmode%
\unskip%
{\dimen1=\baselineskip \divide\dimen1 by 2 %
\raise\dimen1\llap{\tiny -\themycomplaints-}}\fi%
\marginpar{\tiny [\themycomplaints]: ##1}}%
}

\withcomplaints 
\newtheorem{theorem}{Theorem}
\newtheorem{conj}{Conjecture}

\begin{document}

\title {When is a symmetric pin-jointed framework isostatic?}

\author{R. Connelly \\Department of Mathematics, Cornell University,\\
Ithaca, NY 14853, USA\\ P.W. Fowler\\
Department of Chemistry, University of Sheffield,\\
Sheffield S3 7HF, UK\\ S.D. Guest\\ Department of Engineering,
University of Cambridge,
\\Trumpington Street, Cambridge CB2 1PZ, UK\\
B. Schulze and W.J. Whiteley\\ Department of Mathematics and
Statistics, York University,\\ 4700 Keele Street, Toronto, ON M3J
1P3, Canada}

\maketitle

\begin{abstract}
Maxwell's rule from 1864 gives a necessary condition for a
framework to be isostatic in 2D or in 3D. Given a framework
with point group symmetry, group representation theory is
exploited to provide further necessary conditions. This paper
shows how, for an isostatic framework, these conditions imply
very simply stated restrictions on the numbers of those
structural components that are unshifted by the symmetry
operations of the framework. In particular, it turns out that
an isostatic framework in 2D can belong to one of only six
point groups. Some conjectures and initial results are
presented that would give sufficient conditions (in both 2D and
3D) for a framework that is realized generically for a given
symmetry group to be an isostatic framework.
\end{abstract}

\section{Introduction}
This paper deals with \emph{isostatic} frameworks, i.e.,
pin-jointed bar assemblies, commonly referred to in engineering
literature as truss structures, that are both kinematically and
statically determinate. Such systems are minimally infinitesimally
rigid and maximally stress free: they can be termed `just rigid'.
Our ultimate goal is to answer the question posed in the title:
when are symmetric pin-jointed frameworks isostatic?  As a first
step, the present paper provides a series of \emph{necessary}
conditions obeyed by isostatic frameworks that possess symmetry,
and also summarizes conjectures and initial results on
\emph{sufficient} conditions.

Frameworks provide a model that is useful in applications
ranging from civil engineering (Graver, 2001) and the study of
granular materials (Donev et al., 2004) to biochemistry
(Whiteley 2005).  Many of these model frameworks have symmetry.
In applications, both practical and theoretical advantages
accrue when the  framework is isostatic.  In a number of
applications, point symmetry of the framework appears
naturally, and it is therefore of interest to understand the
impact of symmetry on the rigidity of the framework.

Maxwell (1864) formulated a necessary condition for
infinitesimal rigidity, a counting rule for 3D pin-jointed
structures, with an obvious counterpart in 2D; these were later
refined by Calladine (1978).  Laman (1970) provided sufficient
criteria for infinitesimal rigidity in 2D, but there are well
known problems in extending this to 3D (Graver et al., 1993).

The Maxwell counting rule, and its extensions, can be re-cast to
take account of symmetry (Fowler and Guest, 2000) using the
language of point-group representations (see, e.g., Bishop, 1973).
The symmetry-extended Maxwell rule gives additional information
from which it has often been possible to detect and explain
`hidden' mechanisms and states of self-stress in cases where the
standard counting rules give insufficient information (Fowler and
Guest, 2002; 2005, Schulze 2008a).  Similar symmetry extensions
have been derived for other classical counting rules (Ceulemans
and Fowler, 1991; Guest and Fowler, 2005).

In the present paper, we will show that the symmetry-extended
Maxwell rule can be used to provide necessary conditions for a
finite framework possessing symmetry to be stress-free and
infinitesimally rigid, i.e., isostatic.  It turns out that
symmetric isostatic frameworks must obey some simply stated
restrictions on the counts of structural components that are fixed
by various symmetries. For 2D systems, these restrictions imply
that isostatic structures must have symmetries belonging to one of
only six point groups. For 3D systems, all point groups are
possible, as convex triangulated polyhedra (isostatic by the
theorems of Cauchy and Dehn (Cauchy 1813, Dehn 1916)) can be
constructed in all groups (Section~\ref{sec:3Diso}), although
restrictions on the placement of structural components may still
apply.

For simplicity in this presentation, we will restrict our
configurations to realisations in which all joints are
distinct. Thus, if we consider an abstract representation of
the framework as a graph, with vertices corresponding to
joints, and edges corresponding to bars, then we are assuming
that the mapping from the graph to the geometry of the
framework is injective on the vertices. Complications can arise
in the non-injective situation, and will be considered
separately (Schulze, 2008a).

The structure of the paper is as follows:  Maxwell's rule, and its
symmetry extended version, are introduced in
Section~\ref{sec:back}, where a symmetry-extended version of a
necessary condition for a framework to be isostatic is given,
namely the equisymmetry of the representations for mechanisms and
states of self-stress.  In Section~\ref{sec:calc} the calculations
are carried out in 2D, leading to restrictions on the symmetries
and configurations of 2D isostatic frameworks, and in 3D, leading
to restrictions on the placement of structural components with
respect to symmetry elements. In Section~\ref{sec:laman} we
conjecture sufficient conditions for a framework realized
generically for a symmetry group to be isostatic, both in the
plane and in 3D.

\section{Background} \label{sec:back}

\subsection{Scalar counting rule}

Maxwell's rule (Maxwell, 1864) in its modern form (Calladine,
1978),
\begin{equation}
    m - s = 3j - b - 6,
    \label{eq:calladine}
    \end{equation}
expresses a condition for the determinacy of an unsupported,
three-dimensional pin-jointed frame, in terms of counts of
structural components.  In equation~(\ref{eq:calladine}), $b$ is
the number of bars, $j$ is the number of joints, $m$ is the number
of infinitesimal internal mechanisms and $s$ is the number of
states of self-stress.  A statically determinate structure has
$s=0$; a kinematically determinate structure has $m=0$;
isostatic structures have $s=m=0$.

The form of (\ref{eq:calladine}) arises from a comparison of the
dimensions of the underlying vector spaces that are associated
with the equilibrium, or equivalently the compatibility,
relationships for the structure (Pellegrino and Calladine, 1986).

Firstly, the equilibrium relationship can be written as
\begin{equation*}
    \mb{A t} = \mb{f}
    \end{equation*}
where $\mb{A}$ is the \emph{equilibrium} matrix;  $\mb{t}$ is a
vector of internal bar forces (tensions), and lies in a vector
space of dimension $b$;  $\mb{f}$ is an assignment of
externally applied forces, one to each joint,  and, as there
are $3j$ possible force components, $\mb{f}$ lies in a vector
space of dimension $3j$ (this vector space is the tensor
product of a $j$-dimensional vector space resulting from
assigning a scalar to each joint, and a $3$-dimensional vector
space in which a 3D force vector can be defined).
Hence $\mb{A}$ is a $3j \times b$ matrix.

A state of self-stress is a solution to $\mb{A t} = \mb{0}$, i.e., a vector
in the nullspace of $\mb{A}$; if $\mb{A}$ has rank $r$, the dimension
of this nullspace is
\begin{equation}
    s = b - r.
    \label{eq:s}
    \end{equation}

Further, the compatibility relationship can be written as
\begin{equation*}
    \mb{C} \mb{d} = \mb{e}
    \end{equation*}
where $\mb{C}$ is the \emph{compatibility} matrix;  $\mb{e}$ is a
vector of infinitesimal bar extensions, and lies in a vector space
of dimension $b$; $\mb{d}$ is a vector of infinitesimal nodal
displacements, there are $3j$ possible nodal displacements and so
$\mb{d}$ lies in a vector space of dimension $3j$. Hence $\mb{C}$
is a $b \times 3j$ matrix.  In fact, it is straightforward to show
(see e.g., Pellegrino and Calladine, 1986) that $\mb{C}$ is
identical to $\mb{A}^\mathrm{T}$. The matrix $\mb{C}$ is closely
related to the rigidity matrix commonly used in the mathematical
literature: the rigidity matrix is formed by multiplying each row
of $\mb{C}$ by the length of the corresponding bar.  Of particular
relevance here is that fact that the rigidity matrix and $\mb{C}$
have an identical nullspace.

A \emph{mechanism} is a solution to $\mb{A}^\mathrm{T}\mb{d} = 0$,
i.e., a vector in the left-nullspace of $\mb{A}$, and the
dimension of this space is $3j-r$.  However, this space has a
basis comprised of $m$ internal mechanisms and $6$ rigid-body
mechanisms, and hence
\begin{equation}
    m+6 = 3j - r.
    \label{eq:m}
    \end{equation}
Eliminating $r$ from (\ref{eq:s}) and (\ref{eq:m}) recovers Maxwell's equation
(\ref{eq:calladine}).

The above derivation assumes that the system is $3$-dimensional,
but it can be applied to $2$-dimensional frameworks, simply
replacing $3j-6$ by $2j-3$:
\begin{equation}
    m-s = 2j - b - 3.
    \label{eq:peterhouse}
    \end{equation}

\subsection{Symmetry-extended counting rule}
The scalar formula (\ref{eq:calladine}) has been shown (Fowler and
Guest, 2000) to be part of a more general symmetry version of
Maxwell's rule. For a framework with point group symmetry ${\cal
G}$,
 \eqnwithbr{3D:}{
   \Gamma(m) - \Gamma(s)
        = \Gamma(j) \times \Gamma_{xyz} - \Gamma(b) -
        \Gamma_{xyz} - \Gamma_{R_xR_yR_z}
    \label{eq:sm}}
where each $\Gamma$ is known in applied group theory as a
\emph{representation} of ${\cal G}$ (Bishop, 1973), or in
mathematical group theory as a \emph{character} (James and
Liebeck, 2001). For any set of objects $q$,  $\Gamma(q)$ can be
considered as a vector, or ordered set, of the traces of the
transformation matrices $\mb{D}_q(R)$ that describe the
transformation of $q$ under each symmetry operation $R$ that lies
in ${\cal G}$.  In this way, (\ref{eq:sm}) may be considered as a
set of equations, one for each class of symmetry operations in
${\cal G}$. Alternatively, and equivalently, each $\Gamma(q)$ can
be written as the sum of irreducible representations/characters of
${\cal G}$ (Bishop, 1973). In (\ref{eq:sm}) the various sets $q$
are sets of bars $b$, joints $j$, mechanisms $m$ and states of
self-stress $s$; $\Gamma_{xyz}$ and $\Gamma_{R_xR_yR_z}$ are the
translational and rotational representations, respectively.
Calculations using (\ref{eq:sm}) can be completed by standard
manipulations of the character table of the group (Atkins, Child
and Phillips, 1970; Bishop, 1973; Altmann and Herzig, 1994).

The restriction of (\ref{eq:sm}) to $2$-dimensional systems
(assumed to lie in the $xy$-plane) is made by replacing
$\Gamma_{xyz}$ with $\Gamma_{xy}$ and  $\Gamma_{R_xR_yR_z}$ with
$\Gamma_{R_z}$, as appropriate to the reduced set of rigid-body motions.
 \eqnwithbr{2D:}{
    \Gamma(m) - \Gamma(s)
        = \Gamma(j) \times \Gamma_{xy} - \Gamma(b) -
        \Gamma_{xy} - \Gamma_{R_z}
        \label{eq:sm2}}

Examples of the application of (\ref{eq:sm}), (\ref{eq:sm2}),
with detailed working, can be found in Fowler and Guest (2000;
2002; 2005), and further background, giving explicit
transformation matrices, will be found in Kangwai and Guest
(2000).

In the context of the present paper, we are interested in
isostatic systems, which have $m=s=0$, and hence obey the symmetry
condition $\Gamma(m)=\Gamma(s)=0$.  In fact, the symmetry Maxwell
equation (\ref{eq:sm}), (\ref{eq:sm2}) gives the necessary
condition $\Gamma(m)-\Gamma(s)=0$, as it cannot detect the
presence of paired equisymmetric mechanisms and states of self
stress.

The symmetry-extended Maxwell equation corresponds to a set of
$k$ scalar equations, where $k$ is the number of irreducible
representations of ${\cal G}$ (the number of rows in the
character table), or equivalently the number of conjugacy
classes of ${\cal G}$ (the number of columns in the character
table). The former view has been used in previous papers; the
latter will be found useful in the present paper for deriving
restrictions on isostatic frameworks.

\subsection{The need for restrictions}

That existence of symmetry typically imposes restrictions on
isostatic frameworks can be seen from some simple general
considerations. Consider a framework having point-group symmetry
$\cal G$. Suppose that we place all bars and joints freely (so
that no bar or joint is mapped onto itself by any symmetry
operation). Both $b$ and $j$ must then be multiples of $|{\cal
G}|$, the order of the group: $b = \bar{b}|{\cal G}|$, $j =
\bar{j}|{\cal G}|$. Can such a framework be isostatic? Any
isostatic framework obeys the scalar Maxwell rule with $m-s=0$ as
a necessary condition.  In three dimensions, we have $b = 3j-6$,
and hence:
 \eqnlabel{3D:}{
    \bar{b} = 3\bar{j} - \frac{6}{|{\cal G}|}.
    }
In two dimensions, we have $b = 2j-3$, and hence:
 \eqnlabel{2D:}{
    \bar{b} = 2\bar{j} - \frac{3}{|{\cal G}|}.
    }
As $\bar{b}$ and $\bar{j}$ are integers, $|{\cal G}|$ is
restricted to values $1$, $2$, $3$ and $6$ in 3D, and $1$ and
$3$ in 2D. Immediately we have that if the point group order is
not one of these special values, it is impossible to construct
an isostatic framework with all structural components placed
freely: any isostatic framework with $|{\cal G}| \ne 1, 3$ (2D)
or $|{\cal G}| \ne 1, 2, 3, 6$ (3D) must have some components
in {\it special positions} (components that are unshifted by
some symmetry operation).

In the Schoenflies notation (Bishop, 1973), the point groups of
orders 1, 2, 3 and 6 are
\begin{center}\begin{tabular}{ll}
    $|{\cal G}| = 1:$   &   ${\cal C}_1$    \\
    $|{\cal G}| = 2:$   &   ${\cal C}_2$, ${\cal C}_s$, ${\cal C}_i$   \\
    $|{\cal G}| = 3:$   &   ${\cal C}_3$    \\
    $|{\cal G}| = 6:$   &   ${\cal C}_{3h}$, ${\cal C}_{3v}$,
        ${\cal D}_3$, ${\cal S}_6$
    \end{tabular}\end{center}

Further restrictions follow from the symmetry-adapted Maxwell
rules (\ref{eq:sm}), (\ref{eq:sm2}). In a hypothetical
framework where all bars and joints are placed freely, the bar
and joint representations are
\begin{equation*}
    \Gamma(b) = \bar{b}\,\Gamma_\mathrm{reg}; \quad \Gamma(j) =
        \bar{j}\,\Gamma_\mathrm{reg}
    \end{equation*}
where $\Gamma_\mathrm{reg}$ is the \emph{regular representation}
of $\cal G$ with trace $|{\cal G}|$ under the identity operation,
and $0$ under all other operations. The representations
$\Gamma_{xyz}$ and $\Gamma_{xy}$ have trace 3 and 2, respectively,
under the identity operation, and hence equations (\ref{eq:sm})
and (\ref{eq:sm2}) become
 \eqnlabel{3D:}{
     \Gamma(m) - \Gamma(s) = 3\bar{j}\Gamma_\mathrm{reg}
        - \bar{b}\Gamma_\mathrm{reg}
        - \Gamma_{xyz} - \Gamma_{R_xR_yR_z},
    }
 \eqnlabel{2D:}{
    \Gamma(m) - \Gamma(s) = 2\bar{j}\Gamma_\mathrm{reg}
        - \bar{b}\Gamma_\mathrm{reg}
        - \Gamma_{xy} - \Gamma_{R_z},
    }
which can be written as,
 \eqnlabel{3D:}{
     \Gamma(m) - \Gamma(s) = \left(3\bar{j} - \bar{b} -
        \frac{6}{|{\cal G}|}\right)
        \Gamma_\mathrm{reg}
    + \left(\frac{6}{|{\cal G}|}\Gamma_\mathrm{reg} - \Gamma_{xyz} -
        \Gamma_{R_xR_yR_z}\right),
    }
 \eqnlabel{2D:}{
    \Gamma(m) - \Gamma(s) = \left(2\bar{j} - \bar{b} -
        \frac{3}{|{\cal G}|}\right)
        \Gamma_\mathrm{reg}
    + \left(\frac{3}{|{\cal G}|}\Gamma_\mathrm{reg} - \Gamma_{xy} -
        \Gamma_{R_z}\right).
    }
If we have arranged that our hypothetical framework has
satisfied the scalar Maxwell rule, we are left with a `discrepancy
term', given by
 \eqnlabel{3D:}{
    \Gamma(m) - \Gamma(s) = \left(\frac{6}{|{\cal G}|}\Gamma_\mathrm{reg}
        - \Gamma_{xyz} -
        \Gamma_{R_xR_yR_z}\right),
    }
 \eqnlabel{2D:}{
    \Gamma(m) - \Gamma(s) = \left(\frac{3}{|{\cal G}|}\Gamma_\mathrm{reg}
        - \Gamma_{xy} -
        \Gamma_{R_z}\right).
    }
Thus in both 2D and 3D, our hypothetical framework cannot be
isostatic unless the rigid-body motions span a multiple of the
regular representation, when the discrepancy term will
disappear. Within groups of the specified orders, this term
disappears only for: in 3D, $\{ {\cal C}_1, {\cal C}_s, {\cal
C}_i, {\cal C}_3, {\cal C}_{3h}, {\cal C}_{3v}, {\cal S}_6\}$,
and in 2D, $\{ {\cal C}_1, {\cal C}_3\}$. Thus, for example, 3D
frameworks of ${\cal C}_2$ or ${\cal D}_3$ symmetry, with all
structural components shifted by all symmetry operations,
cannot be isostatic, even when they satisfy the scalar Maxwell
count: in both cases, evaluation of the discrepancy term shows
that the hypothetical 3D framework would have a totally
symmetric mechanism and a state of self stress that is
antisymmetric with respect to two-fold rotation.

Frameworks of higher symmetry, such as the icosahedral ($|{\cal
G}|$ = 120 or 60) or cubic groups ($|{\cal G}|$ = 48, 24 or 12)
cannot satisfy even the isostatic scalar Maxwell count without
having structural components in special positions.

Thus, even this simple example shows that for many groups some
restriction on positions of bars and points is imposed by
symmetry, and implies that symmetry adds extra necessary
conditions for frameworks to be isostatic.

\section{Derivation of conditions for isostatic frameworks}
    \label{sec:calc}
In order to apply (\ref{eq:sm}) to any particular framework, we
require, in addition to the standard $\Gamma_{xyz}$ and $\Gamma
_{R_xR_yR_z}$, a knowledge of the bar and joint permutation
representations: $\Gamma(b)$ and $\Gamma(j)$.  In other words,
for each symmetry operation in ${\cal G}$, we need to determine
the numbers of bars and joints that remain unshifted by that
operation. It is necessary to perform this count only once per
conjugacy class.

Setting $\Gamma(m)-\Gamma(s)$ to zero in (\ref{eq:sm}) and
(\ref{eq:sm2}), class by class, will give up to $k$ independent
necessary conditions for the framework to be isostatic.  We
will carry out this procedure once and for all for all point
groups, as there is a very limited set of possible operations
to consider. The two-dimensional and three-dimensional cases
will be considered separately.

\subsection{Two-dimensional isostatic frameworks}

In this section we treat the two-dimensional case:  bars, joints,
and their associated displacements are all confined to the plane.
(Note that frameworks that are isostatic in the plane may have
out-of-plane mechanisms when considered in $3$-space.)  The
relevant symmetry operations are: the identity ($E$), rotation by
$2\pi/n$ about a point ($C_n$), and reflection in a line
($\sigma$). The possible groups are the groups ${\cal C}_n$ and
${\cal C}_{nv}$ for all natural numbers $n$. ${\cal C}_n$ is the
cyclic group generated by $C_n$, and ${\cal C}_{nv}$ is generated
by a $\{ C_n, \sigma\}$ pair.  The group ${\cal C}_{1v}$ is
usually called ${\cal C}_s$.

All two-dimensional cases can be treated in a single calculation,
as shown in Table~\ref{tab:2D}. Each entry in the table is the
trace (character) of the appropriate representation (indicated in
the left column) of the symmetry (indicated in the top line).
Characters are calculated for four operations: we distinguish
$C_2$ from the $C_n$ operation with $n>2$.  Each line in the table
represents a stage in the evaluation of (\ref{eq:sm2}).  Similar
tabular calculations are found in Fowler and Guest (2000) and
subsequent papers.

\begin{figure}
    \begin{center}
        \includegraphics{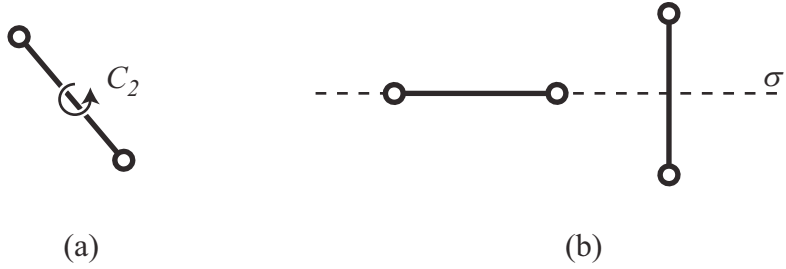}
        \end{center}
    \caption{Possible placement of a bar with respect to a symmetry element
        in two dimensions, such that it is unshifted by the associated symmetry operation:
        (a) $C_{2}$ centre of rotation;  (b) mirror line.}
    \label{fig:2Dbars}
    \end{figure}

To treat all two-dimensional cases in a single calculation, we
need a notation that keeps track of the fate of structural
components under the various operations, which in turn depends on
how the joints and bars are placed with respect to the symmetry
elements. The notation used in Table~\ref{tab:2D} is as follows.
\begin{description}
\item[$j$] is the total number of joints;
\item[$j_c$] is the number of joints lying on the point of
    rotation ($C_{n>2}$ or $C_2$)
    (note that, as we are considering only cases where all joints are distinct,
    $j_{c}= 0$ or $1$);
\item[$j_\sigma$] is the number of joints lying on a given
    mirror line;
\item[$b$] is the total number of bars;
\item[$b_2$] is the number of bars left unshifted by a $C_2$ operation
    (see Figure~\ref{fig:2Dbars}(a) and note that $C_n$ with $n>2$
    shifts all bars);
\item[$b_\sigma$] is the number of bars unshifted by a
    given mirror operation (see Figure~\ref{fig:2Dbars}(b):
    the unshifted bar may lie in, or perpendicular to, the
    mirror line).
\end{description}
Each of the counts refers to a particular symmetry element and
any structural component may therefore contribute to one or
more count, for instance, a joint counted in $j_c$ also
contributes to $j_\sigma$ for each mirror line present.
\begin{table}\begin{center}
    \begin{tabular}{l|c c c c}
                        & $\phantom{-}E$       &  $C_{n > 2}$   &
                            $\phantom{-}C_2$   & $\phantom{-}\sigma$  \\ \hline
    $\phantom{=}\Gamma(j)$& $\phantom{-}j$     & $j_c$             &
                            $\phantom{-}j_c$   & $\phantom{-}j_\sigma$\\
    $\phantom{=}\times \Gamma_{xy}$&   $\phantom{-}2$     & $2\cos\phi$       &
                            $-2$    & $\phantom{-}0$       \\ \hline
    $=\Gamma(j)\times\Gamma_{xy}$
                        & $\phantom{-}2j$      & $2 j_c \cos\phi$  &
                            $-2 j_c$& $\phantom{-}0$       \\
    $\phantom{=}-\Gamma(b)$        & $-b$      & $0$               &
                            $-b_2$  & $-b_\sigma$\\
    $\phantom{=}-(\Gamma_{xy}+\Gamma_{R_z})$
                        & $-3$      & $-2\cos\phi-1$   &
                            $\phantom{-}1$     & $\phantom{-}1$ \\ \hline
    $=\Gamma(m) - \Gamma(s)$
                        & $2j-b-3$  & $2(j_c-1)\cos\phi - 1$ &
                            $-2j_c - b_2 +1$ & $-b_\sigma+1$

    \end{tabular}
    \caption{Calculations of characters for
        representations for the 2D
        symmetry-extended Maxwell equation (\ref{eq:sm2}).}
    \label{tab:2D}
    \end{center}
\end{table}

>From Table~\ref{tab:2D}, the symmetry treatment of the 2D Maxwell
equation reduces to scalar equations of four types.  If $\Gamma(m)
- \Gamma(s) = 0$, then
 \eqnwithbr{$E$:}
    {2j-b=3
    \label{eq:2DE}}
 \eqnwithbr{$C_2$:}
    {2j_c + b_2 = 1
    \label{eq:2DC2}}
 \eqnwithbr{$\sigma$:}
    {b_\sigma = 1
    \label{eq:2Dsigma}}
 \eqnwithbr{$C_{n>2}$:}
    {2(j_c-1)\cos\phi = 1
    \label{eq:2DCn}}
where a given equation applies when the corresponding symmetry
operation is present in ${\cal G}$. Some observations on 2D
isostatic frameworks, arising from this set of equations are:
\begin{enumerate}
\item
Trivially, all 2D frameworks have the identity element and
(\ref{eq:2DE}) simply restates the scalar Maxwell rule
(\ref{eq:peterhouse}) with $m-s=0$.

\item Presence of a $C_2$ element imposes limitations on the
placement of bars and joints.  As both $j_c$ and $b_2$ must be
non-negative integers, (\ref{eq:2DC2}) has the unique solution
$b_2 = 1$, $j_c = 0$.  In other words, an isostatic 2D framework
with a $C_2$ element of symmetry has no joint on the point of
rotation, but exactly one bar centred at that point.

\item
Similarly, presence of a mirror line implies, by (\ref{eq:2Dsigma}), that $b_\sigma = 1$ for that
line, but places no restriction on the number of joints in the same line, and hence allows this bar
to lie either in, or perpendicular to, the mirror.

\item
Deduction of the condition imposed by a rotation of higher order
$C_{n>2}$ proceeds as follows.  Equation (\ref{eq:2DCn}) with $\phi
= 2 \pi/n$ implies
\begin{equation}
    (j_c - 1)\cos\left(\frac{2 \pi}{n}\right) = \frac{1}{2}
    \label{eq:Cn}
    \end{equation}
and as $j_c$ is either $0$ or $1$, this implies that $j_c =0$ and $n=3$. Thus, a 2D isostatic
framework cannot have a $C_n$ rotational element with $n>3$, and when either a $C_2$ or $C_3$
rotational element is present, no joint may lie at the centre of rotation.
\end{enumerate}

In summary, a 2D isostatic framework may have only symmetry
operations drawn from the list $\{E,\, C_2,\, C_3,\, \sigma\}$, and
hence the possible symmetry groups ${\cal G}$ are 6 in number:
${\cal C}_1$, ${\cal C}_2$, ${\cal C}_3$, ${\cal C}_s$, ${\cal
C}_{2v}$, ${\cal C}_{3v}$.  Group by group, the conditions necessary
for a 2D framework to be isostatic are then as follows.
\begin{description}
\item[${\cal C}_1$:] $b=2j-3$.
\item[${\cal C}_2$:] $b=2j-3$ with $b_2 = 1$ and $j_c = 0$, and
    as all other bars and joints occur in pairs, $j$ is even and
    $b$ is odd.
\item[${\cal C}_3$:] $b=2j-3$ with $j_c = 0$, and hence all joints
    and bars occur in sets of $3$.
\item[${\cal C}_s$:] $b=2j-3$ with $b_\sigma = 1$ and all other bars
    occurring in pairs.  Symmetry does not restrict $j_\sigma$.
\item[${\cal C}_{2v}$:] $b=2j-3$ with $j_c = 0$ and $b_2 = b_\sigma = 1$.
    A central bar lies in one of the two mirror lines, and
    perpendicular to the other.  Any additional bars must
    lie in the general position, and hence occur
    in sets of $4$, with joints
    in sets of $2$ and $4$.  Hence  $b$ is odd and $j$ is even.
\item[${\cal C}_{3v}$:] $b=2j-3$ with $j_c = 0$ and $b_\sigma = 1$
    for each of the three mirror lines.
\end{description}
We consider whether these condition are also sufficient in
Section~\ref{sec:2Dcond}.

\begin{figure}
    \begin{center}
        \includegraphics{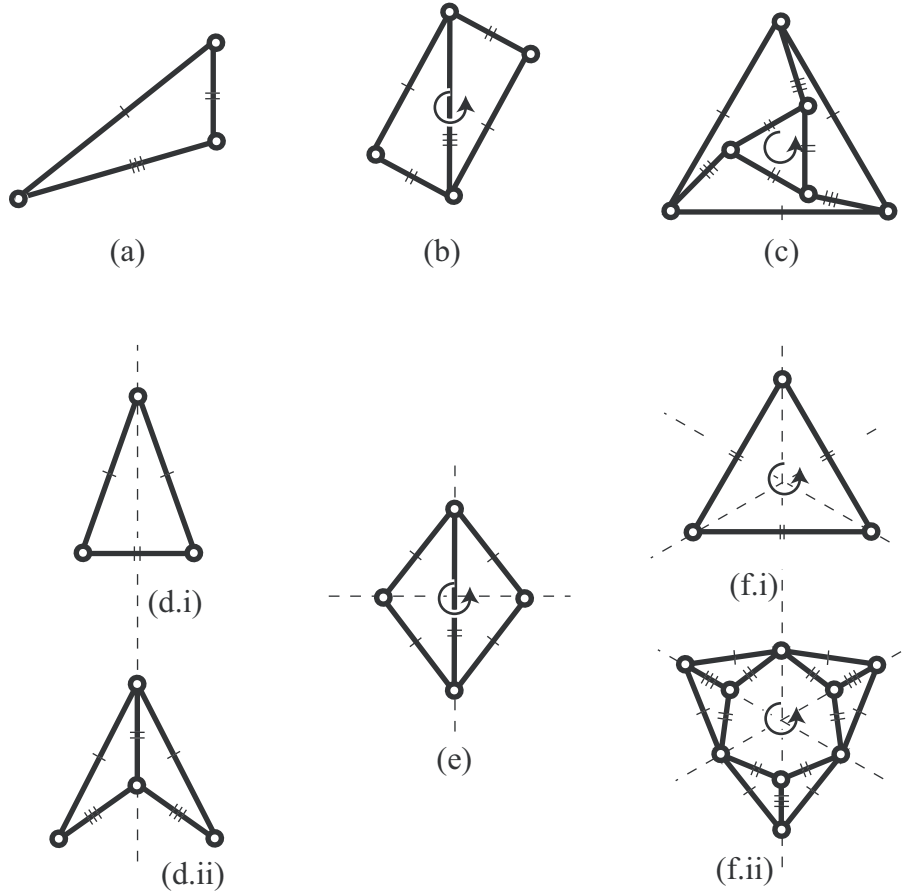}
        \end{center}
    \caption{Examples, for each of the possible groups,
    of small 2D isostatic
    frameworks, with bars which are equivalent under symmetry marked with the same symbol:
    (a) ${\cal C}_1$; (b) ${\cal C}_2$;
    (c) ${\cal C}_3$; (d) ${\cal C}_s \equiv {\cal C}_{1v}$;
    (e) ${\cal C}_{2v}$; (f) ${\cal C}_{3v}$.
    Mirror lines are shown dashed, and rotation axes are indicated
    by a circular arrow.  For each of
    ${\cal C}_s$ and ${\cal C}_{3v}$, two examples are given:
    (i) where each mirror has a bar centered at, and perpendicular to, the mirror line;
    (ii) where a bar lies in each mirror line.   For ${\cal C}_{2v}$, the bar lying
    at the centre must lie in one mirror line, and perpendicular to the
    other.}
    \label{fig:examples}
    \end{figure}

Figure~\ref{fig:examples} gives examples of small 2D isostatic
frameworks for each of the possible groups, including cases
where bars lie in, and perpendicular to, mirror lines.

\subsection{Three-dimensional isostatic frameworks}
\label{sec:3Diso}

The families of possible point groups of 3D objects are: the
icosahedral ${\cal I}$, ${\cal I}_{h}$; the cubic ${\cal T}$,
${\cal T}_{h}$, ${\cal T}_{d}$, ${\cal O}$, ${\cal O}_{h}$; the
axial ${\cal C}_n$, ${\cal C}_{nh}$, ${\cal C}_{nv}$; the dihedral
${\cal D}_n$, ${\cal D}_{nh}$, ${\cal D}_{nd}$; the cyclic ${\cal
S}_{2n}$; and the trivial ${\cal C}_s$, ${\cal C}_i$, ${\cal
C}_{1}$ (Atkins et al., 1970).  The relevant symmetry operations
are: proper rotation  by $2\pi/n$ about an axis, $C_n$, and
improper rotation, $S_{n}$ ($C_n$ followed by reflection in a
plane perpendicular to the axis). By convention, the identity $E
\equiv C_1$, inversion $i \equiv S_2$, and reflections $\sigma
\equiv S_1$ are treated separately.

The calculation is shown in Table~\ref{tab:3D}.  Characters are
calculated for six operations.  For proper rotations, we
distinguish $E$ and $C_2$ from the $C_n$ operations with $n>2$.
For improper rotations, we distinguish $\sigma$ and $i$ from
the $S_{n>2}$ operations. We exclude from consideration the
degenerate case of a single bar, and assume that the total
number of joints is greater than three.

\begin{figure}
    \begin{center}
        \includegraphics{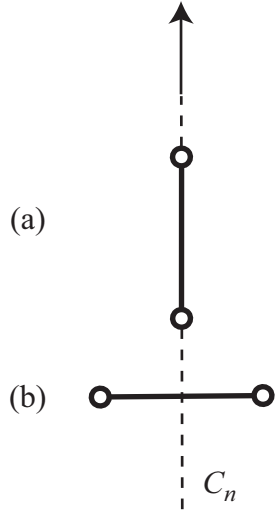}
        \end{center}
    \caption{Possible placement of a bar unshifted by a proper
    rotation about an axis:
    (a) for any $C_{n \geq 2}$;  (b) for $C_{2}$ alone.}
    \label{fig:3Daxis}
    \end{figure}

\begin{figure}
    \begin{center}
        \includegraphics{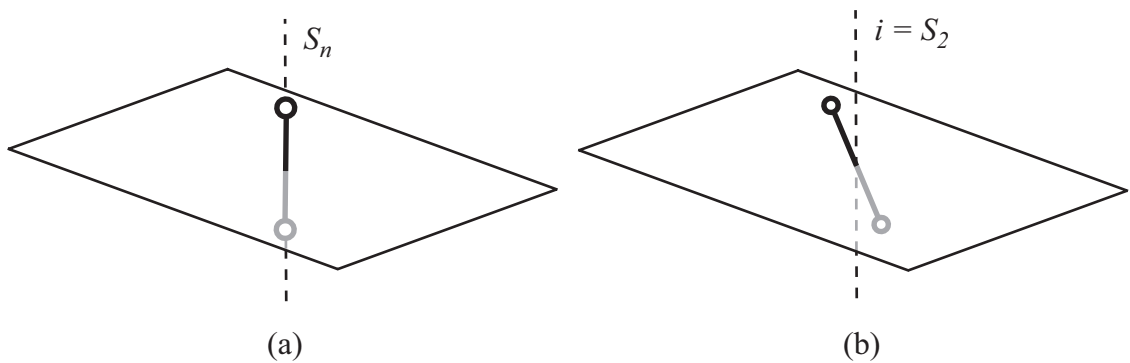}
        \end{center}
    \caption{Possible placement of a bar unshifted by an improper
    rotation about an axis:
    (a) for any $S_{n \geq 2}$;  (b) for $i=S_{2}$.}
    \label{fig:3DS}
    \end{figure}

\begin{figure}
    \begin{center}
        \includegraphics{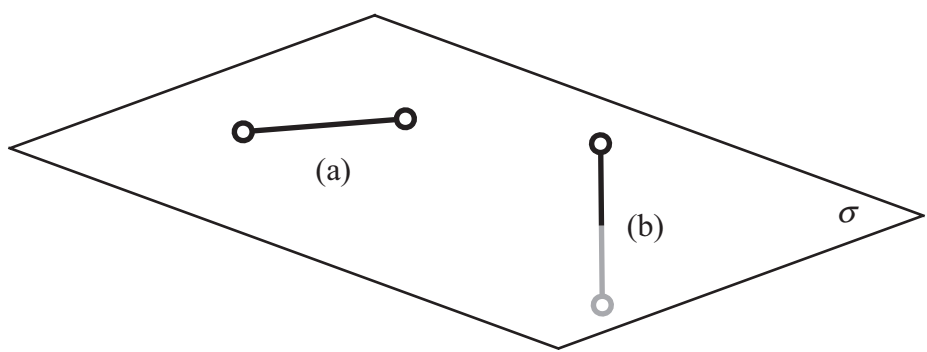}
        \end{center}
    \caption{Possible placement of a bar unshifted by a reflection
    in  a plane: (a) lying in the plane; (b) lying perpendicular to
    the plane.}
    \label{fig:3Dreflect}
    \end{figure}

The notation used in Table~\ref{tab:3D} is
\begin{description}
\item[$j$] is the total number of joints;
\item[$j_n$] is the number of joints lying on the $C_n$ axis;
\item[$j_c$] is the number of joints ($0$ or {1}) lying on the unique central
    point (if any).  Such joints are unshifted by all operations;
\item[$j_\sigma$] is the number of joints lying on a given $\sigma$
    mirror-plane;
\item[$b$] is the total number of bars;
\item[$b_n$] is the number of bars unshifted by a $C_{n>2}$
    rotation:
    note that each such bar must lie along the axis of the rotation
    (see Figure~\ref{fig:3Daxis}(a));
\item[$b_{nc}$] is the number of bars unshifted by the improper
    rotation $S_{n>2}$: note that such bars must lie along the
    axis of the rotation, and be centered on the central point of the group
    (see Figure~\ref{fig:3DS}(a));
\item[$b_c$] is the number of bars unshifted by the inversion $i$:
    note that the centre of the bar must lie at the central point
    of the group, but no particular orientation is implied
    (see Figure~\ref{fig:3DS}(b));
\item[$b_2$] is the number of bars unshifted by the $C_{2}$
    rotation:
    such bars must lie either along, or perpendicular to and centered on,
    the axis (see Figure~\ref{fig:3Daxis}(a) and (b));
\item[$b_\sigma$] is the number of bars unshifted by a given $\sigma$
    mirror operation (see Figure~\ref{fig:3Dreflect}(a) and (b)).
\end{description}
Again, each of the counts refers to a particular symmetry element, and so, for instance the joint
counted in $j_c$ also contributes to $j$, $j_n$ and $j_\sigma$.
\begin{sidewaystable}\begin{center}
    \begin{tabular}{l|c c c c c c c c}
                        & $\phantom{-}E$       &  $C_{n \neq 2}(\phi)$ &
                          $\phantom{-}C_2$     & $\phantom{-}\sigma$            &
                          $\phantom{-}i$                   &
                          $S_{n \neq 2}(\phi)$  \\ \hline
    $\phantom{=}\Gamma(j)$& $\phantom{-}j$     & $\phantom{-}j_n$                 &
                          $\phantom{-}j_2$     & $\phantom{-}j_\sigma$            &
                          $\phantom{-}j_c$                 &
                          $\phantom{-}j_c$                 \\
    $\phantom{=}\times \Gamma_{xyz}$&$\phantom{-}3$       & $2\cos\phi+1$         &
                          $-1$      & $\phantom{-}1$                   &
                          $-3$                  &
                          $2\cos\phi-1$         \\ \hline
    $=\Gamma(j)\times \Gamma_{xyz}$
                          & $\phantom{-}3j$    & $(2\cos\phi+1)j_n$    &
                          $-j_2$    & $\phantom{-}j_\sigma$            &
                          $-3j_c$               &
                          $(2\cos\phi-1)j_c$    \\
    $\phantom{=}-\Gamma(b)$        & $-b$      & $-b_n$                &
                          $-b_2$    & $-b_\sigma$           &
                          $-b_c$                &
                          $-b_{nc}$             \\
    $\phantom{=}-(\Gamma_{xyz}+\Gamma_{R_xR_yR_z})$
                        & $-6$      & $-4\cos\phi-2$        &
                          $\phantom{-}2$       & $\phantom{-}0$                   &
                          $\phantom{-}0$                   &
                          $\phantom{-}0$                   \\ \hline
    $=\Gamma(m) - \Gamma(s)$
                        & $3j-b-6$  & $(2\cos\phi+1)(j_n-2) - b_n$ &
                          $-j_2 - b_2 +2$ & $j_\sigma-b_\sigma$ &
                          $-3j_c-b_c$ &
                          $(2\cos\phi-1)j_c - b_{nc}$
    \end{tabular}
    \caption{Calculations of characters for
        representations in the 3D
        symmetry-extended Maxwell equation (\ref{eq:sm}).}
    \label{tab:3D}
    \end{center}
\end{sidewaystable}

{\samepage
 From Table~\ref{tab:3D}, the symmetry treatment of the
3D Maxwell equation reduces to scalar equations of six types. If
$\Gamma(m) - \Gamma(s) = 0$, then
 \eqnwithbr{$E$:}
    {3j-b=6
    \label{eq:3DE}}
 \eqnwithbr{$\sigma$:}
    {b_\sigma = j_\sigma
    \label{eq:3Dv}}
 \eqnwithbr{$i$:}
    {3j_c + b_c = 0
    \label{eq:3Di}}
 \eqnwithbr{$C_2$:}
    {j_2 + b_2 = 2
    \label{eq:3DC2}}
 \eqnwithbr{$C_{n>2}$:}
    {(j_n-2)(2\cos\phi+1) = b_n
    \label{eq:3DCn}}
 \eqnwithbr{$S_{n>2}$:}
    {j_c(2\cos\phi-1) = b_{nc}
    \label{eq:3DSn}}
where a given equation applies when the corresponding symmetry
operation is present in ${\cal G}$. }

Some observations on 3D isostatic frameworks, arising from the
above, are:
\begin{enumerate}
\item From (\ref{eq:3DE}), the framework must satisfy the
scalar Maxwell rule (\ref{eq:calladine}) with $m-s=0$.
\item From (\ref{eq:3Dv}), each mirror that is present
    contains the same number of joints as bars that are
    unshifted under reflection in that mirror. \item From
    (\ref{eq:3Di}), a centro-symmetric framework has
    neither a joint nor a bar centered at the inversion
    centre.
    \item For a $C_2$ axis, (\ref{eq:3DC2}) has solutions
\begin{equation*}
    (j_2, b_2) = (2,0), (1,1), (0,2).
    \end{equation*}
The count $b_2$ refers to both bars that lie along, and
those that lie perpendicular to, the axis.  However, if a
bar were to lie along the $C_2$ axis, it would contribute
$1$ to $b_2$ and $2$ to $j_2$ thus generating a
contradiction of (\ref{eq:3DC2}), so that in fact all bars
included in $b_2$ must lie perpendicular to the axis. \item
Equation (\ref{eq:3DCn}) can be written, with $\phi =
2\pi/n$, as
\begin{equation*}
    (j_n-2)\left(2\cos\left(\frac{2 \pi}{n}\right)+1\right) = b_n
    \end{equation*}
with $n>2$.  The non-negative integer solution $j_n = 2$, $b_n =
0$, is possible for all $n$.  For $n>2$ the factor
$(2\cos(2\pi/n)+1)$ is rational at $n=3, 4, 6$, but generates a
further distinct solution only for $n=3$:
\begin{description}
    \item[$n=3$]
    \begin{equation*}
        0(j_3 - 2) = b_3
        \end{equation*}
    and so here $b_3 = 0$, but $j_3$ is unrestricted.
    \item[$n=4$]
    \begin{equation*}
        j_4 - 2 = b_4
        \end{equation*}
    $C_4$ implies $C_4^2 = C_2$ about the same axis, and hence
    $b_4 = 0$, and $j_4 = j_2 = 2$.
    \item[$n=6$]
    \begin{equation*}
        2(j_6 - 2) = b_6
        \end{equation*}
    $C_6$ implies $C_6^3 = C_2$ and $C_6^2 = C_3$ about the
    same axis, and hence $b_6 = b_3 = 0$, and $j_6 = j_3 = j_2 = 2$.
    \end{description}
Thus $b_n$ is $0$ for any $n>2$, and only in the case $n=3$ may
$j_n$ depart from $2$.
\item
Likewise, equation (\ref{eq:3DSn}) can be written, with $\phi =
2\pi/n$, as
\begin{equation*}
    \left(2\cos\left(\frac{2 \pi}{n}\right)-1\right) j_c = b_{nc}
    \end{equation*}
with $n>2$.  The integer solution $j_c = 0$, $b_{nc} = 0$, is
possible for all $n$.  For $n>2$ the factor $(2\cos(2\pi/n)-1)$ is
rational at $n=3, 4, 6$, but generates no further solutions:
\begin{description}
    \item[$n=3$]
    \begin{equation*}
        -2j_c = b_{3c}
        \end{equation*}
    and so $j_c = b_{3c} = 0$.
    \item[$n=4$]
    \begin{equation*}
        -j_c = b_{4c}
        \end{equation*}
    and so $j_c = b_{4c} = 0$.
    \item[$n=6$]
    \begin{equation*}
        0j_c = b_{6c}
        \end{equation*}
    and $b_{6c} = 0$
    but $S_6$ implies $S_6^3 = i$ and hence also $j_c = 0$.
    \end{description}

\item
For a framework with icosahedral (${\cal I}$ or
${\cal I}_h$) symmetry, the requirement that $j_5 = 2$ for
each 5-fold axis implies that the framework must include a
single orbit of $12$ vertices that are the vertices of an
icosahedron. Similarly, for a framework with a $\cal{O}$ or
${ \cal O}_h$ symmetry, the requirement that $j_4 = 2$
implies that the framework must include a single orbit of
$6$ vertices that are the vertices of an octahedron.
\end{enumerate}

\begin{figure}
    \hspace*{\fill} \includegraphics[height=5cm]{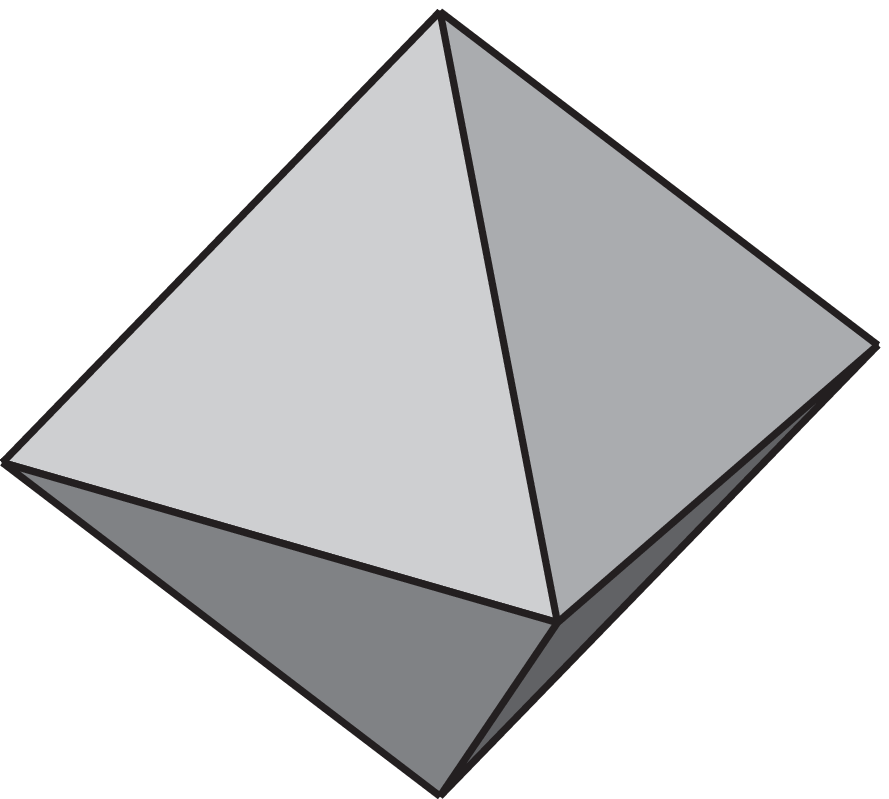}
    \hspace*{\fill} \includegraphics[height=5cm]{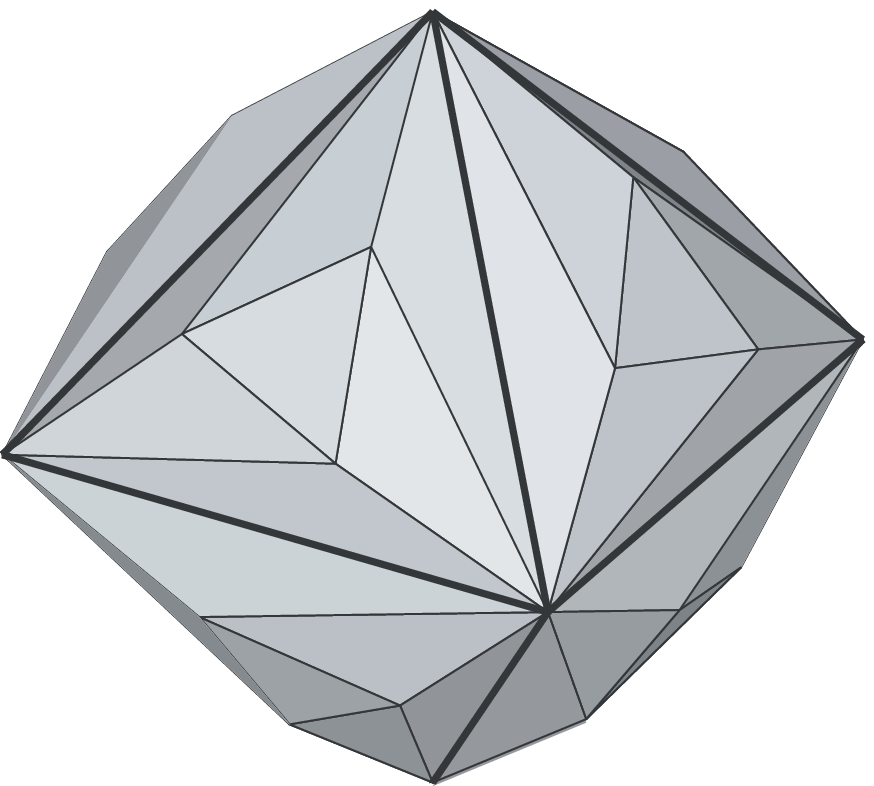} \hspace*{\fill}\\[12pt]
    \hspace*{\fill} (a) \hspace*{5.8cm} (b) \hspace*{\fill}
    \caption{A regular octahedron (a), and a convex polyhedron (b) generated
        by adding a twisted octahedron to every face of the original
        octahedron.  The polyhedron in (b) has the rotation but
        not the reflection symmetries of the polyhedron in (a).
        If a framework is constructed from either polyhedron by placing bars
        along edges, and joints at vertices, the framework will be isostatic.}
     \label{fig:twistcap}
   \end{figure}

In contrast to the 2D case, in 3D the symmetry conditions do
not exclude any point group.  For example, a fully triangulated
convex polyhedron, isostatic by the Theorem of Cauchy and Dehn
(Cauchy 1813; Dehn 1916) can be constructed to realize any 3D
point group. Beginning with the regular triangulated polyhedra
(the tetrahedron, octahedron, icosahedron), infinite families
of isostatic frameworks can be constructed by expansions of
these polyhedra using operations of truncation and capping. For
example, to generate isostatic frameworks with only the
rotational symmetries of a given triangulated polyhedron, we
can `cap' each face with a twisted octahedron, consistent with
the rotational symmetries of the underlying polyhedron: the
resultant polyhedron will be an isostatic framework with the
rotational symmetries of the underlying polyhedron, but none of
the reflection symmetries. An example of the capping of a
regular octahedron is shown in Figure~\ref{fig:twistcap}.
Similar techniques can be applied to create polyhedra for any
of the point groups.

\begin{figure}
    \hspace*{\fill} {\includegraphics[height=5cm]{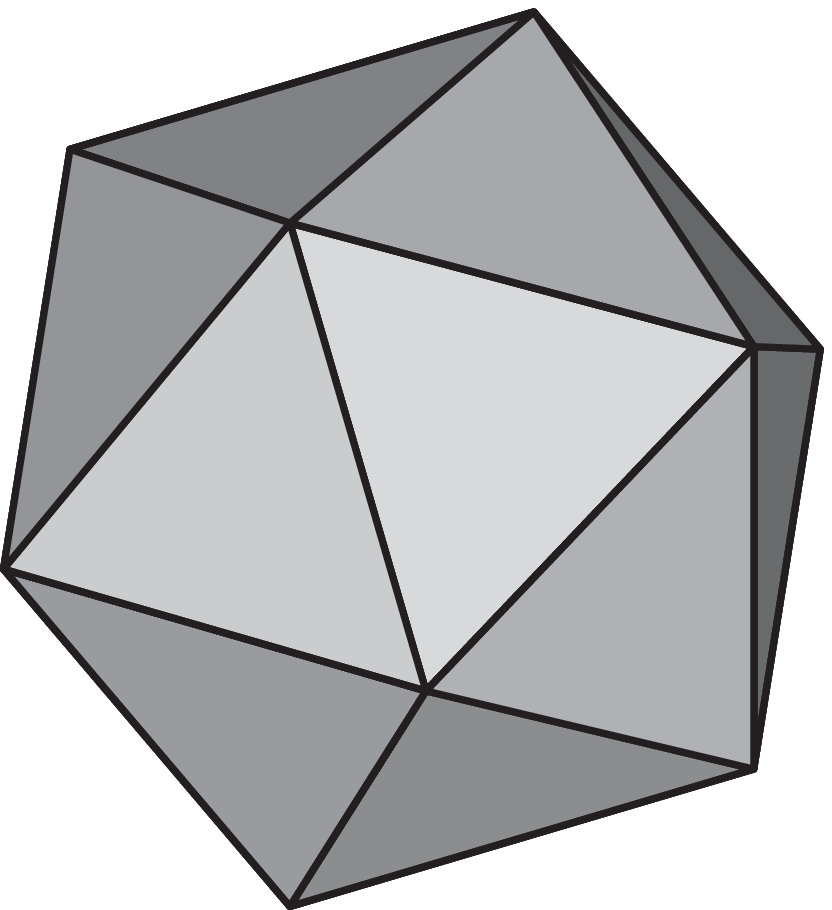}}
    \hspace*{\fill} {\includegraphics[height=5cm]{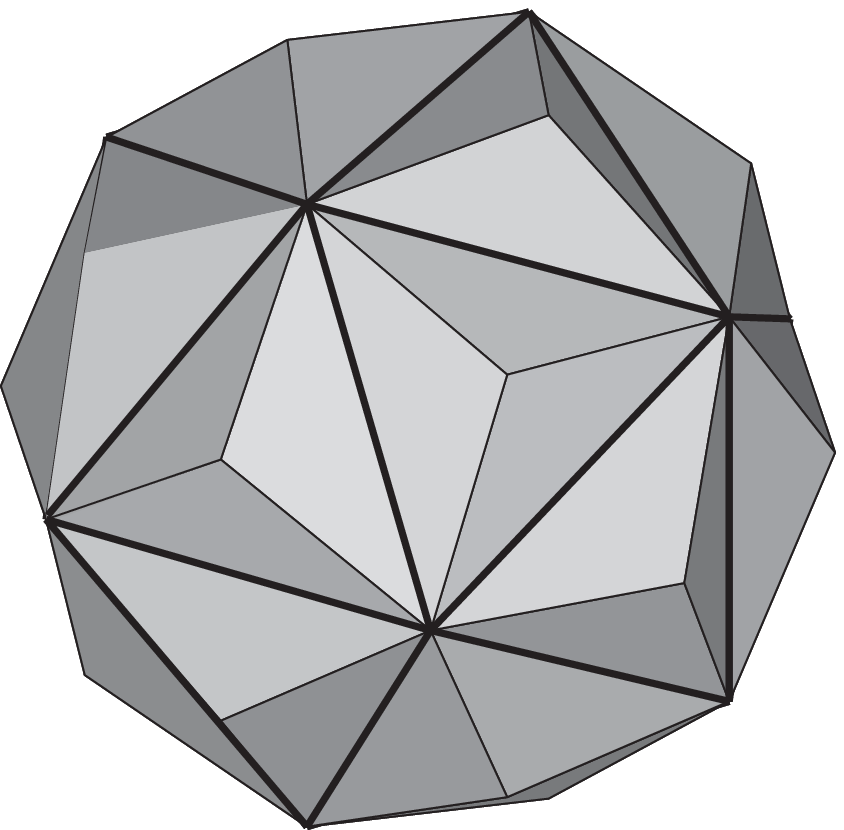}} \hspace*{\fill} \\[12pt]
    \hspace*{\fill} (a) \hspace*{5.8cm} (b) \hspace*{\fill}
    \caption{An icosahedron (a), and the second stellation of the
    icosahedron (b).  If a framework is constructed from either
    polyhedron by placing bars along edges, and joints at vertices,
    the framework will be isostatic.  The framework (b) could be
    constructed from the framework (a) by `capping' each face
    of the original icosahedron preserving the $C_{3v}$ site symmetry.}
    \label{fig:ico}
    \end{figure}

\begin{figure}
    \begin{center}
        \includegraphics{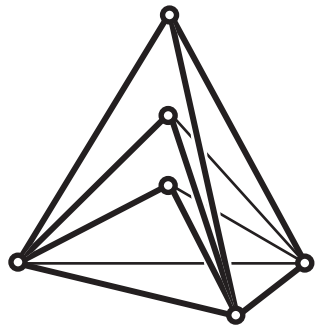}
        \end{center}
    \caption{A series of `hats' added symmetrically along a 3-fold
    axis of an isostatic framework leaves the structure isostatic.}
    \label{fig:cap}
    \end{figure}

One interesting possibility arises from consideration of groups
that contain $C_3$ axes. Equation (\ref{eq:3DCn}) allows an
unlimited number of joints, though not bars, along a 3-fold
symmetry axis. Thus, starting with an isostatic framework,
joints may be added symmetrically along the 3-fold axes.  To
preserve the Maxwell count, each additional joint is
accompanied by $3$ new bars. Thus, for instance, we can `cap'
every face of an icosahedron to give the compound
icosahedron-plus-dodecahedron (the second stellation of the
icosahedron), as illustrated in Figure~\ref{fig:ico}, and this
process can be continued ad infinitum adding a pile of `hats'
consisting of a new joint, linked to all three joints of an
original icosahedral face (Figure~\ref{fig:cap}). Similar
constructions starting from cubic and trigonally symmetric
isostatic frameworks can be envisaged. Addition of a single
`hat' to a triangle of a framework is one of the Hennenberg
moves (Tay \& Whiteley 1985): changes that can be made to an
isostatic framework without introducing extra mechanisms or
states of self stress.

\section{Sufficient Conditions for Isostatic Realisations.}
  \label{sec:laman}

\subsection{Conditions for two-dimensional isostatic frameworks}
\label{sec:2Dcond}

For a framework with point-group symmetry ${\cal G}$ the
previous section has provided some necessary conditions for the
realization to be isostatic. These conditions included some
over-all counts on bars and joints, along with sub-counts on
special classes of bars and joints (bars on mirrors or
perpendicular to mirrors, bars centered on the axis of
rotation, joints on the centre of rotation etc.). Here,
assuming that the framework is realized with the joints in a
configuration as generic as possible (subject to the symmetry
conditions), we investigate whether these conditions are
sufficient to guarantee that the framework is isostatic.

The simplest case is the identity group (${\cal C}_1$). For
this basic situation,  the key result is Laman's Theorem.  In
the following, we take $G=\{J,B\}$ to define the connectivity
of the framework, where $J$ is the set of $j$ joints and $B$
the set of $b$ bars, and we take $p$ to define the positions of
all of the joints in 2D.
\begin{theorem}\label{Laman} {(Laman, 1970)} For a generic configuration
in 2D, $p$, the framework $G(p)$ is isostatic if and only if
$G=\{J,B\}$ satisfies the conditions:
\begin{enumerate}
\item $b = 2j -3$;
\item for any non-empty set of bars $B^{*}$, which contacts just the joints in $J^{*}$,
with $|B^{*}|=b^{*}$ and
$|J^{*}|=j^{*}$,  $b^{*} \leq  2j^{*} -3$.
\end{enumerate}
\end{theorem}

Our goal is to extend these results to other symmetry groups.
With the appropriate definition of `generic' for symmetry
groups (Schulze 2008a), we can anticipate that the necessary
conditions identified in the previous sections for the
corresponding group plus the Laman condition identified in
Theorem \ref{Laman}, which considers subgraphs that are not
necessarily symmetric, will be sufficient. For three of the
plane symmetry groups, this has been confirmed. We use the
previous notation for the point groups and the identification
of special bars and joints, and describe a configuration as
`generic with symmetry group ${\cal G}$' if, apart from
conditions imposed by symmetry, the points are in a generic
position (the constraints imposed by the local site symmetry
may remove 0,1 or 2 of the two basic freedoms of the point).

\begin{theorem}  (Schulze 2008b) If $p$ is a plane configuration generic
with symmetry group ${\cal G}$, and $G(p)$ is a framework realized with these symmetries,
then the following  necessary conditions are also sufficient for $G(p)$ to be isostatic:
\vskip .1in  \noindent  {$b = 2j -3$} and for any non-empty set
of bars $B^{*}$,   $b^{*} \leq 2j^{*} -3$ and

\begin{enumerate}
 \item   for $C_{s}$:   $b_{\sigma} = 1$;
\item   for $C_{2}$:   $b_{2} = 1$, $j_{c}=0$
\item   for $C_{3}$:   $ j_{c}=0$
   \end{enumerate}

\end{theorem}

For the remaining groups, we have a conjecture.

\begin{conj}  If $p$ is a plane configuration generic with symmetry
group ${\cal G}$, and $G(p)$ is a framework realized with these symmetries,
then the following  necessary conditions are also sufficient for $G(p)$ to be isostatic:
\vskip .1in  \noindent  {$b = 2j -3$} and for any non-empty set
of bars $B^{*}$,   $b^{*} \leq 2j^{*} -3$ and

\begin{enumerate}
\item   for $C_{2v}$:  $b_{2} =1$ and $b_{{\sigma}}= 1$ for
    each mirror
\item   for $C_{3v}$:   $ j_{c}=0$ and $b_{\sigma} = 1$ for
    each mirror.
   \end{enumerate}

\end{conj}

An immediate consequence of this theorem (and the conjecture)
is that there is (would be) a polynomial time algorithm to
determine whether a given framework in generic position modulo
the symmetry group ${\cal G}$ is isostatic. Although the Laman
condition of Theorem \ref{Laman} involves an exponential number
of subgraphs of $G$, there are several algorithms that
determine whether it holds in $cjb$ steps where $c$ is a
constant.  The pebble game (Hendrickson and Jacobs, 1997) is an
example. The additional conditions for being isostatic with the
symmetry group ${\cal G}$ trivially can be verified in constant
time.

\subsection{Conditions for isostatic 3D frameworks}
In 3D, there is no known counting characterization of generically
isostatic frameworks, although we have the necessary conditions:
$3j - b -6 =0$ and $3j' - b' -6  \geq  0$ for all subgraphs with
$j' \geq  3$ (Graver 2001). There are, however a number of
constructions for graphs which are known to be generically
isostatic in 3D (see, for example, Tay and Whiteley 1985, Whiteley
1991). If we assume that we start with such a graph $G$, then it
is natural to ask whether the additional necessary conditions for
a realization $G(p)$ that is generic with point group symmetry
${\cal G}$ to be isostatic are also sufficient. In contrast to the
plane case, where we only needed to state these conditions once,
for the entire graph, in 3D  for all subgraphs $G'$ of $G$ whose
realizations $G'(p)$ are symmetric with a subgroup ${\cal G}'$ of
${\cal G}$, with the full count $3j'-b'-6=0$, we need to assert
the conditions corresponding to the symmetry operations in ${\cal
G}'$ as well. These conditions are clearly necessary, and for all
reflections, half-turns, and $6$-fold rotations in ${\cal G}'$,
they do not follow from the global conditions on the entire graph
(as they would in the plane). See Schulze (2008c) for details.
\\\indent All of the above conditions
combined, however, are still not sufficient for a $3$-dimensional
framework $G(p)$ which is generic with point group symmetry ${\cal
G}$ to be isostatic, because even if $G(p)$ satisfies all of these
conditions, the symmetry imposed by ${\cal G}$ may force parts of
$G(p)$ to be `flattened' so that a self-stress of $G(p)$ is
created. For more details on how `flatness' caused by symmetry
gives rise to additional necessary conditions for $3$-dimensional
frameworks to be isostatic, we refer the reader to Schulze,
Watson, and Whiteley (2008).

\section*{References}
Altmann, S.L. and Herzig, P., 1994. Point-Group Theory Tables.
Clarendon Press, Oxford. \\[6pt]
Atkins, P.W., Child, M.S. and Phillips, C.S.G., 1970. Tables for
Group Theory. Oxford University Press, Oxford. \\[6pt]
Bishop, D.M., 1973.  Group Theory and Chemistry. Clarendon Press, Oxford. \\[6pt]
Calladine, C.R., 1978. Buckminster Fuller's `Tensegrity'
structures and Clerk Maxwell's rules for the construction of stiff
frames.
International Journal of Solids and Structures 14, 161--172.\\[6pt]
Cauchy, A.L., 1813. Recherche sur les poly\`{e}dres --- premier
m\'{e}moire, Journal de l'Ecole Polytechnique 9, 66--86.\\[6pt]
Ceulemans, A. and Fowler, P.W., 1991. Extension of Euler's theorem
to the
symmetry properties of polyhedra. Nature 353, 52--54. \\[6pt]
Dehn, M., 1916. \"{U}ber die Starreit konvexer Polyeder,
Mathematische Annalen 77, 466--473.\\[6pt]
Donev A., Torquato, S., Stillinger, F. H.  and Connelly, R., 2004.
Jamming in hard sphere and disk packings.
Journal of Applied Physics 95, 989--999.\\[6pt]
Fowler, P.W. and Guest, S.D., 2000. A symmetry extension of
Maxwell's rule for rigidity of frames.  International Journal of
Solids and Structures 37, 1793--1804. \\[6pt]
Fowler, P.W. and Guest, S.D., 2002.  Symmetry and states of self
stress in triangulated toroidal frames.  International Journal of
Solids and Structures 39, 4385--4393.  \\[6pt]
Fowler, P.W. and Guest, S.D., 2005.  A symmetry analysis of
mechanisms in rotating rings of tetrahedra.  Proceedings of the
Royal Society: Mathematical, Physical \& Engineering Sciences.
461(2058), 1829-1846. \\[6pt]
Guest, S.D. and Fowler, P.W. 2005. A symmetry-extended mobility
rule. Mechanism and Machine Theory. 40, 1002-1014. \\[6pt]
Graver, J.E, Servatius, B.,  and Servatius, H., 1993.
Combinatorial Rigidity.
Graduate Studies in Mathematics, AMS, Providence. \\[6pt]%
Graver, J.E, 2001. Counting on Frameworks: Mathematics to Aid
the Design of Rigid Structures.
The Mathematical Association of America, Washington, DC. \\[6pt]
Hendrickson, B. and Jacobs, D., 1997.
An algorithm for two-dimensional rigidity percolation: The pebble game.
 Journal of Computational Physics, 137, 346--365\\[6pt]
James, G. and Liebeck, M., 2001.  Representations and Characters
of Groups, 2nd edition.  Cambridge University Press. \\[6pt]
Kangwai, R.D. and Guest, S.D., 2000. Symmetry-adapted equilibrium
matrices.  International Journal of Solids and Structures 37, 1525--1548. \\[6pt]
Laman, G., 1970. On graphs and rigidity of plane skeletal structures.
Journal of Engineering Mathematics 4, 331--340. \\[6pt]
Maxwell, J.C., 1864. On the calculation of the equilibrium and
stiffness of frames, Philosophical Magazine 27, 294--299. Also:
Collected papers, XXVI. Cambridge University Press,
1890.\\[6pt]
Pellegrino, S. and Calladine, C.R., 1986. Matrix analysis of
statically and kinematically indeterminate structures.
International Journal of Solids and Structures 22, 409--428.
\\[6pt]
Schulze, B., 2008a.  Injective and non-injective realizations
with symmetry. Preprint,
{\tt arXiv:0808.1761} \\[6pt]
Schulze, B., 2008b.  Symmetrized Laman's Theorems, in
preparation,
York University, Toronto, Canada.\\[6pt]
Schulze, B., 2008c.  Combinatorial and geometric rigidity with
symmetry constraints, PhD thesis in
preparation,
York University, Toronto, Canada.\\[6pt]
Schulze, B., Watson, A., and Whiteley, W., 2008.  Symmetry, flatness, and
necessary conditions for independence, in
preparation.\\[6pt]
Tay, T-S. and Whiteley, W., 1985. Generating Isostatic
Frameworks. Structural Topology 11, 20--69.
\\[6pt]
Whiteley, W., 1991 Vertex splitting in isostatic frameworks,
Structural Topology 16, 23--30.
\\[6pt]
Whiteley, W., 2005. Counting out the flexibility of proteins.
Physical Biology 2, 116--126.

\end{document}